\documentclass[8pt,a4paper]{article}
\usepackage[dvips]{graphicx}
\usepackage{latexsym}
\usepackage{amssymb}
\usepackage{amsmath}
\usepackage{amscd}
\newtheorem{Th}{Theorem}[section]
\newtheorem{Co}[Th]{Corollary}
\newtheorem{Lem}[Th]{Lemma}

\newtheorem{Pro}[Th]{Proposition}
\newcommand{\demo}{\par\noindent{\it Proof. \/}\ }
\newcommand{\enD}{\hfill $\Box$\vspace{3truemm} \par}
\newcommand{\bx}{\mbox{\boldmath $x$}}
\newcommand{\bX}{\mbox{\boldmath $X$}}
\newcommand{\be}{\mbox{\boldmath $e$}}
\newcommand{\ba}{\mbox{\boldmath $a$}}
\newcommand{\bb}{\mbox{\boldmath $b$}}
\newcommand{\bv}{\mbox{\boldmath $v$}}
\newcommand{\by}{\mbox{\boldmath $y$}}
\newcommand{\bo}{\mbox{\boldmath $0$}}

\newcommand{\bw}{\mbox{\boldmath $w$}}

\newcommand{\blambda}{\mbox{\boldmath $\lambda$}}
\newcommand{\bn}{\mbox{\boldmath $n$}}
\newcommand{\R}{{\mathbb R}}
\newcommand{\lon}{\longrightarrow}
\begin{document}
\title{\bf The $S_t^1\times S_s^1$-valued lightcone Gauss map
of a Lorentzian surface  in semi-Euclidean 4-space
 }
\author{\normalsize \bf  Donghe Pei, Lingling Kong, Jianguo Sun and Qi Wang}
\date{\today}
\maketitle
\begin{abstract}
We define the notions of $S_t^1\times S_s^1$-valued lightcone Gauss
maps, lightcone pedal surface  and  Lorentzian lightcone height
function of Lorentzian surface in semi-Euclidean 4-space and
established the relationships between singularities of these objects
and geometric invariants of the surface  as applications of standard
techniques of singularity theory for the Lorentzian lightcone height
function.
\end{abstract}
\renewcommand{\thefootnote}{\fnsymbol{footnote}}
 \footnote[0]{Work partially supported
by NSF of China No.10471020 and NCET of China No.05-0319. }
\section{Introduction}
In \cite{Izu3,Izu4}, S.Izumiya et al studied singularities of
lightcone Gauss maps and lightlike hypersurfaces of spacelike
surface in Minkowski 4-space,
 and established the relationships between
such singularities  and geometric invariants of these surfaces under
the action of Lorentz group. Our aim in this paper is to develop the
analogous study  for Lorentzian surface  in semi-Euclidean 4-space
$\R_2^4$. To do this we need to develop first the local differential
geometry  of Lorentzian surface  in semi-Euclidean 4-space $\R_2^4$
in  a similar way than the classically done surfaces in Euclidean
4-space \cite{Oneil}. As it was to be expected, the situation
presents certain peculiarities when compared with the Euclidean
case.  For instance, in our case it is always possible to choose two
lightlike normal directions along the Lorentzian surface  a frame of
its normal bundle. By using this, we define a Lorentzian invariant
${\mathcal K}_l(1,\pm 1)$ and call it the {\it lightlike
Gauss-Kronecker curvature} of the Lorentzian surface. We introduce
the notion of lightcone height function and use it to show that the
$S_t^1\times S_s^1$-valued lightcone Gauss map has a singular point
if and only if the lightlike Gauss-Kronecker curvature vanishes at
such point. Moreover, we show that the $S_t^1\times S_s^1$-valued
lightcone Gauss map is a constant map if and only if the Lorentzian
surface is contained in a lightlike hyperplane, so we can view the
singularities of the $S_t^1\times S_s^1$-valued lightcone Gauss map
as an estimate of the contacts of the surface  with lightlike
hyperplanes.

\par
We shall assume throughout the whole paper that all the
maps and manifolds are $C^{\infty}$ unless the contrary is explicitly stated.

\par
Let $\Bbb R^4=\{(x_1,x_2,x_3,x_4)|x_i \in \R\ (i=1,2,3,4)\ \}$ be a
4-dimensional vector space. For any vectors $\bx=(x_1,x_2,x_3,x_4)$
and $\by =(y_1,y_2,y_3,y_4)$ in $\Bbb R^4,$ the {\it pseudo scalar
product \/} of $\bx$ and $\by$  is defined to be
$\langle\bx,\by\rangle =-x_1 y_1-x_2 y_2 +x_3y_3+x_4 y_4$. We call
$(\Bbb R^4, \langle ,\rangle )$ a {\it semi-Euclidean 4-space\/} and
write $\Bbb R^4_2$ instead of $(\Bbb R^4,\langle ,\rangle )$.
\par
We say that a vector $\bx$ in $\Bbb R^4_2\setminus\{\bo \}$ is {\it
spacelike\/}, {\it lightlike\/} or {\it timelike\/} if
$\langle\bx,\bx\rangle >0,=0$ or $<0$ respectively. The norm of the
vector $\bx\in \Bbb R_2^4$ is defined by $\Vert \bx \Vert =\sqrt{
|\langle\bx,\bx\rangle |}.$ For a lightlike vector $\bn\in \Bbb
R^4_2$ and a real number $c$, we define the {\it lightlike
hyperplane\/} with pseudo normal $\bn$ by
$$LHP(\bn,c)=\{\bx \in \Bbb R^4_2 |\langle\bx,\bn\rangle
=c\}.$$
\par
Let $\bX: U\to \Bbb R^4_2$ an immersion, where $U\subset \Bbb R^2$
is an open subset. We denote that $M=\bX (U) $ and identify $M$ and
$U$ by the immersion $\bX$.
\par
We say that $M$ is a {\it Lorentzian surface \/} if the tangent
space $T_p  M$ of $M$ is a Lorentzian surface  for any point $p\in
M$. In this case, the normal space $N_p M$ is a Lorentzian plane.
Let $\{\be_3 (x,y),\be_4 (x,y);p=(x,y)\}$ be an pseudo-orthonormal
frame of the tangent space $T_p M$ and $\{\be_1 (x,y),\be_2
(x,y);p=(x,y)\}$ a pseudo-orthonormal frame of $N_p M$, where,
$\be_1(p),\be_3(p)$ are unit timelike vectors and $\be_2,\be_4$ are
unit spacelike vectors.
\par
We shall now establish the fundamental formula for a Lorentzian
2-space in $\Bbb R_2^4$ by means of similar notions to those of [9].

\par
We  can write $d\bX=\sum\limits_{i=1}^4 \omega _i \be_i$ and
$d\be_i=\sum\limits_{j=1}^4 \omega _{ij} \be_j;\ i=1,2,3,4.$ where
$\omega _i$ and $\omega _{ij}$ are $1$-forms given by $\omega _i=
\delta (\be_i)\langle d\bX,\be_i\rangle$ and $\omega _{ij}=\delta
(\be_j)\langle d\be_i,\be_j\rangle ,$
\par
with \ \ \ \ \ \ \ \ \ \ \ \ \ $\delta ({\mathbf{e}}_{i})=<{\mathbf{e}}_{i},{\mathbf{e}}_{i}>=\left\{%
\begin{array}{ll}
    1, & \hbox{$i=2,4,$} ;\\
    -1, & \hbox{$i=1,3.$}. \\
\end{array}%
\right.    $
\par
We have the Codazzi type equations:
\begin{eqnarray}
\left\{%
\begin{array}{ll}
    d\omega_{i}=\sum _{j=1}^4\delta({\mathbf{e}}_{i})\delta({\mathbf{e}}_{j})\omega_{ij}\wedge \omega_{j};\\
    d\omega_{ij}=\sum _{k=1}^4\omega_{ik} \wedge \omega_{kj}, \\
\end{array}%
\right.
\end{eqnarray}
\par
where $d$ is exterior derivative.

\par

Since $\langle \be_i,\be_j\rangle =\delta _{ij}\delta (\be_j)$
(where $\delta _{ij}$ is Kronecker's delta),
we get
\begin{eqnarray}
\omega_{ij}=-\delta (\be_i)\delta (\be_j)\omega_{ji}.
\end{eqnarray}
In particular, $\omega_{ii}=0;\ i=1,2,3,4.$ It follows from the fact
$\langle d\bX,\be_1\rangle =\langle d\bX,\be_2\rangle =0$ that
\begin{eqnarray}
\omega _1=\omega _2 =0.
\end{eqnarray}
\par
Therefore we have
\begin{eqnarray}
\left\{
  \begin{array}{ll}
    0=d\omega _1 = \delta(e_j)\sum\limits_{j=1}^4
\omega_{1j}\wedge\omega_{j}= \delta(e_j)\sum\limits_{j=3}^4
\omega_{1j}\wedge\omega_{j}=
-\omega_{13}\wedge\omega_{3}+\omega_{14}\wedge\omega_{4}; \\
    0=d\omega _{2} =\delta(e_j)\sum\limits_{j=1}^4
\omega_{2j}\wedge\omega_{j}=\delta(e_j)\sum\limits_{j=3}^4
\omega_{2j}\wedge\omega_{j}=
-\omega_{23}\wedge\omega_{3}+\omega_{24}\wedge\omega_{4}.
  \end{array}
\right.
\end{eqnarray}

\par
By Cartan's lemma, we can write
\begin{eqnarray}
 \left\{
  \begin{array}{ll}
    \omega_{13}=a\omega_{3}+b\omega_{4};\ \ \ \ \omega_{14}=-b\omega_{3}-c\omega_{4}; \\
    \omega_{23}=\overline{a}\omega_{3}+\overline{b}\omega_{4};\ \ \ \ \omega_{24}=-\overline{b}\omega_{3}-\overline{c}\omega_{4}.
  \end{array}
\right.
\end{eqnarray}
for appropriate functions $a,b,c,\bar a,\bar b,\bar c.$

\par
Since $\langle d\bX,\be_1\rangle =\langle d\bX,\be_2\rangle =0$,
\begin{eqnarray*}
\langle d^2\bX,\be_1\rangle &=&-\langle d\bX,d\be_1\rangle \\
&=& -\langle\sum\limits_{i=1}^4\omega _i\be_i,
\sum\limits_{j=1}^4\omega _{1j}\be_j\rangle =
-\langle\sum\limits_{i=3}^4\omega _i\be_i,
\sum\limits_{j=2}^4\omega _{1j}\be_j\rangle\\
&=& -(-\omega _3\omega _{13}+\omega _4\omega _{14})
\\
 &=&
\omega _3(a\omega _{3}+b\omega _4)+\omega _{4}(b\omega _{3}+c\omega
_4)
\\
 &=&a\omega _{3}^2+2b\omega _3\omega _4+c\omega _4^2.
\end{eqnarray*}
\par
As the same $\langle d^2\bX,\be_2\rangle=\overline{a}\omega
_{3}^2+2\overline{b}\omega _3\omega _4+\overline{c}\omega _4^2.$
\par
Then we have a vector-valued quadratic form
$$ -\langle d^2\bX,\be_1\rangle \be_1+\langle d^2\bX,e_2\rangle \be_2
= -(a\omega _{3}^2+c\omega _4^2+2b\omega _3\omega _4)\be_1 +(\bar
a\omega _{3}^2+\bar c\omega _4^2+2\bar b\omega _3\omega _4)\be_2,$$
which is called the second fundamental form of the Lorentz surface.

By using equations (2) and  a straight forward calculation leads us
to the following equations:
$$
d\left(%
\begin{array}{c}
  \mathbf{e_1+e_2} \\
 \mathbf{ e_1-e_2} \\
  \mathbf{e_3} \\
  \mathbf{e_4}\\
\end{array}%
\right)=\left(%
\begin{array}{cccc}
  0 & \omega_{12} & \omega_{13}+\omega_{23} & \omega_{14}+\omega_{24} \\
  -\omega_{12} & 0 & \omega_{13}-\omega_{23}& \omega_{14}-\omega_{24} \\
  -(\omega_{13}+\omega_{23})/2 & (\omega_{23}-\omega_{13})/2 & 0 & \omega_{34} \\
  (\omega_{14}+\omega_{24})/2 & (\omega_{14}-\omega_{24})/2  & \omega_{34} & 0 \\
\end{array}%
\right)
\left(%
\begin{array}{c}
 \mathbf{ e_1-e_2} \\
  \mathbf{e_1+e_2} \\
  \mathbf{e_3} \\
 \mathbf{ e_4}\\
\end{array}%
\right)\\
$$
\par

On the other hand, we define
$$LC_p=\{\bx=(x_1,x_2,x_3,x_4)\in\Bbb R_2^4|\ -(x_1-p_1)^2-(x_2-p_2)^2+
(x_3-p_3)^2+(x_4-p_4)^2=0\}$$ and $$S_t^1\times
S_s^1=\{\bx=(x_1,x_2,x_3,x_4) \in LC_0\ |\ x_1^2+x_2^2=1, \ x_1\geq
0,\ x_2\geq 0 \},$$ where $p=(p_1,p_2,p_3,p_4)\in \Bbb R_2^4$,
$S_t^1$ denotes the {\it timelike circle \/} and $S_s^1$ denotes the
{\it spacelike circle  \/}. We call $LC_p^*=LC_p\setminus\{p\}$ a
{\it lightcone \/} at the vertex $p$. Given any lightlike vector
$\bx= (x_1,x_2,x_3,x_4)$, we have $\widetilde{\bx}=
(\frac{x_1}{\sqrt{x_1^2+x_2^2}},\frac{x_2}{\sqrt{x_1^2+x_2^2}},
\frac{x_3}{\sqrt{x_1^2+x_2^2}},\frac{x_4}{\sqrt{x_1^2+x_2^2}} )\in
S_t^1\times S_s^1$.
\par
Let $\be_1=(a_1,a_2,a_3,a_4)$, $\be_2=(b_1,b_2,b_3,b_4)$, and $\xi
^{\pm}=\sqrt{(a_1-b_1)^2\pm (a_2-b_2)^2}$, then we have the
following fundamental formula:
$$d\left(%
\begin{array}{c}
  \widetilde{\mathbf{e_1+e_2}} \\
  \widetilde{\mathbf{e_1-e_2}} \\
  \mathbf{e_3} \\
  \mathbf{e_4}\\
\end{array}%
\right)=\left(%
\begin{array}{cccc}
  0 & -\omega_{12}-\frac{d\xi^{+}}{\xi^{+}} & \frac{\omega_{13}+\omega_{23}}{\xi^{+}} & \frac{\omega_{14}+\omega_{24}}{\xi^{+}} \\
  -\omega_{12}-\frac{d\xi^{-}}{\xi^{-}} & 0 & \frac{\omega_{13}-\omega_{23}}{\xi^{-}}& \frac{\omega_{14}-\omega_{24}}{\xi^{-}} \\
  -(\omega_{13}+\omega_{23})/2 & (\omega_{23}-\omega_{13})/2 & 0 & \omega_{34} \\
  (\omega_{14}+\omega_{24})/2 & (\omega_{14}-\omega_{24})/2  & \omega_{34} & 0 \\
\end{array}%
\right)
\left(%
\begin{array}{c}
 \widetilde{ \mathbf{e_1-e_2}} \\
  \widetilde{\mathbf{e_1+e_2}} \\
  \mathbf{e_3} \\
  \mathbf{e_4}\\
\end{array}%
\right)$$\\

\par
Given $\bv=x\be_1+y\be_2\in N_pM,$ we have
$d\bv=dx\be_1+xd\be_1+dy\be_2+yd\be_2$, and then

$$\langle d\bv,\be_3\rangle\wedge \langle d\bv,\be_4\rangle ={\mathcal K}_l(x,y)
\omega _3\wedge\omega_4,$$ where the function ${\mathcal K}_l$ as
follows:
\begin{eqnarray*}
{\mathcal K}_l(x,y)&=&(ax+\bar{a}y)(cx+\bar{c}y) -(bx+\bar{b}y)^2.
\end{eqnarray*}

\par
On the other hand, we define two maps
$$LG_M^\pm :U\longrightarrow S_t^1\times S_s^1$$
by $LG_M^\pm (x,y)=\widetilde{\be_1\pm \be_2}(x,y).$ Each one of
these maps shall be called {\it $S_t^1\times S_s^1$-valued lightcone
Gauss map\/} of $\bX(U) =M.$
\par
Now we introduce the notion of Lorentzian lightcone height functions
on the Lorentzian surface  in $\R_2^4$ which is useful for the study
of singularities of $S_t^1\times S_s^1$-valued lightcone Gauss maps.
\par
For a Lorentzian surface $M(=\bX(U))\in\R_2^4$,  we now define a
function
$$H:U\times S_t^1\times S_s^1\longrightarrow \Bbb R$$
by $H((x,y),\blambda)=\langle \bX(x,y),\blambda\rangle$, where
$\blambda=(\cos\theta,\sin\theta,\lambda _3,\lambda _4)\in
S_t^1\times S_s^1$. We call $H$ the {\it Lorentzian lightcone height
function\/} on the surface  $M$. We denote that $h_{\lambda
_0}(x,y)= H(x,y,\blambda _0 )$, for any fixed ${\blambda _0}\in
S_t^1\times S_s^1$. Then we have the following proposition.
\par
\begin{Pro}\rm{
Let $M$ be a Lorentzian surface  in $\R_2^4$ and $H:U\times
S_t^1\times S_s^1\longrightarrow \Bbb R$ a Lorentzian lightcone
height function.
Then we have the following assertions:
\par\noindent
{\rm (1)} $(\partial h_\lambda/\partial x)(p_0)= (\partial
h_\lambda/\partial y)(p_0)=0$  if and only if
\par\noindent
$$\blambda = \mu (\be_1\pm\be_2)(p_0)=\widetilde{\be_1
\pm\be_2}(p_0),$$ where $\be_1(p_0)=(a_1,a_2,a_3,a_4)$,
$\be_2(p_0)=(b_1,b_2,b_3,b_4)$ and $\mu =\frac {1}{\sqrt{(a_1\pm
b_1)^2+(a_2\pm b_2)^2}};$ for any point\linebreak
$p_0)4=(x_0,y_0)\in M,$
\par\noindent
{\rm (2)} $(\partial h_\lambda/\partial x)(p_0)= (\partial
h_\lambda/\partial y)(p_0)= {\rm det}{\mathcal
H}(h_{\lambda})(p_0)=0 $
 if and only if
 $$\blambda =\widetilde{\be_1 \pm\be_2}(p_0)\ and
\ {\mathcal K}_l(1,\pm 1)(p_0)=0.$$ Here, ${\rm det}{\mathcal
H}(h_{\lambda})(x,y) $ is the determinant of the Hessian matrix of
$h_{\lambda}$ at $(x,y)$.}
\end{Pro}
\demo By a straight forward calculation, $(\partial
h_\lambda/\partial x)(p_0)= (\partial h_\lambda/\partial y)(p_0)=0$
if and only if
$$
\langle \bX_x,\blambda\rangle (p_0)=\langle \bX_y,\blambda\rangle
(p_0)= 0.
$$
It is equivalent to the condition that $\blambda \in N_{p_0}M$ and
$\blambda\in S_t^1\times S_s^1.$ This means that $\blambda = \mu
(\be_1\pm\be_2)=\widetilde{\be_1\pm\be_2}.$
\par
On the other hand, we now choose local coordinates such that $\bX$
is given by the Monge form $\bX(x,y)=(f_1(x,y),x,f_2(x,y),y)$ and
$\be_1(p_0)=(1,0,0,0)$ and $\be_2(p_0)=(0,0,1,0).$ Since
$${\rm det}{\mathcal H}(h_{\lambda})(x,y) =\vmatrix
\langle \bX_{xx},\blambda\rangle&\langle \bX_{xy},\blambda\rangle \\
\langle \bX_{xy},\blambda\rangle&\langle \bX_{yy},\blambda\rangle
\endvmatrix=0$$ and $\blambda (p_0)=(1,0,\pm1,0) ,$
we have
$$
\vmatrix \langle (f_{1_{xx}},0,f_{2_{xx}},0),\blambda (p_0) \rangle&
\langle (f_{1_{xy}},0,f_{2_{xy}},0),\blambda (p_0)
\rangle\\
\langle (f_{1_{xy}},0,f_{2_{xy}},0),\blambda (p_0) \rangle& \langle
(f_{1_{yy}},0,f_{2_{yy}},0),\blambda (p_0) \rangle
\endvmatrix $$
$$=\vmatrix
a\pm \bar a&b\pm \bar b \\
b\pm \bar b&c\pm \bar c
\endvmatrix=0.
$$
This is equivalent to the condition that ${\mathcal K}_l(1,\pm
1)(x,y)=0$ and $\blambda (p_0)=\widetilde{e_1\pm e_2} .$ \enD
\par

As a corollary of the above proposition, we have the following
theorem.

\begin{Th}\rm{
Under the same assumption as the assumption of the above
proposition, the following conditions are equivalent:
\par\noindent
{\rm (1)} $p\in M$ is a degenerate singular point of Lorentzian
lightcone height function $h_{\lambda}.$
\par\noindent
{\rm (2)}
 There is $\blambda\in S_t^1\times S_s^1$ such that
$(p,\blambda )$ is a singular point of the $S_t^1\times
S_s^1$-valued lightcone Gauss map $LG_M^{\pm}.$
\par\noindent
{\rm (3)} ${\mathcal K}_l(1,\pm 1)(p)=0.$}
\end{Th}
\demo We denote that
$$
\Sigma (H)=\left\{(p,\blambda)\in U\times S_t^1\times S_s^1\ |\
\frac{\partial h_\lambda}{\partial x}(p)= \frac{\partial h_\lambda
}{\partial y}(p)= 0\ \right\}.
$$
By above proposition, (1), we have
$$
\Sigma (H)=\{(p,\blambda )\in U\times S_t^1\times S_s^1\ |\ \blambda
=\widetilde{\be_1\pm \be_2}(p)\ \}.
$$
We now consider the canonical projection $\pi :U\times S_t^1\times
S_s^1\longrightarrow S_t^1\times S_s^1,$ then $\pi |\Sigma (H)$ can
be identified to the $S_t^1\times S_s^1$-valued lightcone Gauss map
$LG^{\pm}_M.$ Under this identification, we can show that the
condition (1) is equivalent to the condition (2).
\par
Above proposition, (2) means  that the condition (2) is equivalent
to the condition (3).
\par
\begin{Th} \rm{Let $M$ be a Lorentzian surface  in $\R_2^4$.
\par\noindent
{\rm (1)} The $S_t^1\times S_s^1$-valued lightcone Gauss maps $
LG_M^{+}$ $($respectively, $LG_M^-)$ is constant if and only if
there exists a unique lightlike hyperplane $LHP(\bv^+,c^+)$
$($respectively, $LHP(\bv^-,c^-))$ such that $M\subset
LHP(\bv^+,c^+)$ $($respectively, $M\subset LHP(\bv^-,c^-))$, where
$\bv^{\pm}=\widetilde{\be_1 \pm\be_2}(x,y)$ and $\langle \bX(x,y),
\bv^{\pm}\rangle =c^{\pm}$ for any $(x,y)\in M$.
\par\noindent
{\rm (2)} Both of the $S_t^1\times S_s^1$-valued lightcone Gauss
maps $ LG_M^{+}$ and $ LG_M^{-}$ are constant if and only if $M$ is
a Lorentzian 2-plane. In this case, the intersection of lightlike
hyperplanes
$$
LHP(\widetilde{\be_1 +\be_2},c^+)\cap LHP(\widetilde{\be_1
-\be_2},c^-)
$$ is the Lorentzian 2-plane $M.$}
\end{Th}
\demo (1) For convenience, we consider the case when
$LG_M^+(x,y)=\widetilde{\be_1 +\be_2}(x,y)$ is constant, so that we
have
$$
d\langle \bX,\widetilde{\be_1 +\be_2}\rangle = \langle
d\bX,\widetilde{\be_1 +\be_2}\rangle + \langle
\bX,d(\widetilde{\be_1 +\be_2})\rangle =0.
$$
Therefore, $\langle \bX,\widetilde{\be_1 +\be_2}\rangle \equiv c^+.$
This means that $M=\bX (U)\subset LHP(\bv^+,c+),$ where
$\bv^+=\widetilde{\be_1 +\be_2}(x,y).$ For the converse assertion,
suppose that there exists a lightlike vector $\bv$ and a real number
$c$ such that $\bX (U)=M\subset LHP(\bv,c).$ Since $\langle
\bX(x,y),\bv\rangle =c,$ we have $d\langle \bX(x,y),\bv\rangle =0.$
This means that $\bv$ is a lightlike normal vector of $M.$ Thus we
have $\widetilde{\bv}=\widetilde{\be_1 \pm\be_2}(x,y).$ This
completes the proof of the assertion (1).
\par
Since $\bv^+\notin LHP(\bv^-,c^-)$ and $\bv^-\notin LHP(\bv^+,c^+),$
$LHP(\bv^-,c^-)$ and $LHP(\bv^+,c^+)$ intersect transversally. By
the assertion (1), both of the $S_t^1\times S_s^1$-valued lightcone
Gauss maps $ LG_M^{+}$ and $ LG_M^{-}$ are constant if and only if
$M\subset LHP(\bv^+,c^+)\cap LHP(\bv^-,c^-).$ Here, the intersection
is a Lorentzian 2-plane. Thus we have the assertion (2). \enD

\par
We say that a point $p_0= (x_0,y_0)$ is a {\it Lorentzian lightlike
parabolic point} of $M$ if ${\mathcal K}_l(1,1)(p_0)=0$ or
${\mathcal K}_l(1,-1)(p_0)=0.$

\section{The lightcone pedal surface  }

\par
In this section we consider a singular hyperplane in the lightcone
$LC_0$ associated to $M$ whose singularities correspond to
singularities of the $S_t^1\times S_s^2$-valued lightcone Gauss map
of $M.$ We now define a family of functions
$$
\widetilde{H}:U\times LC_0\longrightarrow \Bbb R
$$
by
$$
\widetilde{H}((x,y),\bv )=\langle \bX(x,y),\widetilde{\bv}\rangle
-\sqrt{v_1^2+ v_2^2},
$$
where $\bv=(v_1,v_2,v_3,v_4).$ We call $\widetilde{H}$ the {\it
extended Lorentzian lightcone height function \/} of $M=\bX(U).$ As
a corollary of above proposition , we have the following
proposition.
\begin{Pro}\rm{
Let $M$ be a Lorentzian surface  in $\R_2^4$ and
$\widetilde{H}:M\times LC_0\longrightarrow \Bbb R$ the extended
Lorentzian lightcone height function of $M.$ For $p_0=(x_0,y_0)$ and
$\bv_0\in LC_0,$ we have the following:
\par\noindent
{\rm (1)} $\widetilde{H}(p_0,\bv_0)=(\partial \widetilde{H}/\partial
x)(p_0,\bv_0) =(\partial \widetilde{H}/\partial y)(p_0,\bv_0)= 0$ if
and only if
$$
\widetilde{\bv}_0 =\widetilde{\be_1 \pm\be_2}(p_0)\
{\rm and}\ \sqrt{v_1^2+ v_2^2}=\langle \bX(p_0), \widetilde{\be_1 \pm\be_2}(p_0)\rangle.$$
\par\noindent
{\rm (2)}
$$\widetilde{H}(p_0,\bv _0)=\frac{\partial \widetilde{H}}{\partial x}(p_0,\bv_0)=
\frac{\partial \widetilde{H}}{\partial y}(p_0,\bv_0)= {\text
det}{\mathcal H}(\widetilde{h}_{v})(p_0)=0
$$
 if and only if
 $$
\widetilde{\bv}_0 =\widetilde{\be_1 \pm\be_2}(p_0), \ \sqrt{v_1^2+
v_2^2}=\langle \bX(p_0), \widetilde{\be_1 \pm\be_2}(p_0)\rangle \
and \ {\mathcal K}_l(1,\pm 1)(p_0)=0.$$ Here, for a fixed
$\mathbf{\bv}\in LC_0,
\widetilde{H}((x,y),\bv)=\widetilde{h}_{v}(x,y).$}
\end{Pro}

\par
The assertion of proposition 2.1 means that the discriminant set of the
extended Lorentzian
lightcone height function
$\widetilde{H}$ is given by
$$
{\mathcal D}_{\widetilde{H}}=\left\{\bv\ |\ \bv =\langle \bX(x,y),
\widetilde{\be_1\pm\be_2}(x,y)\rangle
(\widetilde{\be_1\pm\be_2})(x,y)\
 \text{for some}\ (x,y)\in U\ \right\}.
$$
Therefore we now define a pair of singular surface  in $LC_0$ by
$$
LP^\pm _M (p)= LP^\pm _M (x,y)=\langle \bX(x,y),
\widetilde{\be_1\pm\be_2}(x,y)\rangle
(\widetilde{\be_1\pm\be_2})(x,y).
$$
We call each $LP^\pm $ the {\it lightcone pedal surface  \/} of
$\bX(U)=M.$ A singularity of the lightcone pedal surface exactly
corresponds to a singularity of the $S_t^1\times S_s^2$-valued
lightcone Gauss map.
\par
We define a pair of hyperplane $LH^{\pm}_{M}:M\times
\mathbb{R}\rightarrow \mathbb{R}_{2}^{4}$~~by
$$LH^{\pm}_{M}(p,u)=LH^{\pm}_{M}(x,y,u)=X(x,y)+u(\widetilde{e_{1}+e_{2}})(x,y),
$$ where $p=X(x,y)$ we call $LH^{\pm}_{M}$ the lightlike hyperplane along M.
\par
We now explain the reason why such a correspondence exists from the
view point of Symplectic and Contact geometry. We consider a point
$\bv=(v_1,v_2,v_3,v_4)\in LC_0,$ then we have a relation
$v_1=\sqrt{-v_2^2+v_3^2+v_4^2}.$ We adopt the coordinate
$(v_2,v_3,v_4)$ of the manifold $LC_0.$ We now consider the
projective cotangent bundle $\pi :PT^*(LC_0)\longrightarrow LC_0$
with the canonical contact structure. We review geometric properties
of this space. Consider the tangent bundle $ \tau
:TPT^*(LC_0)\rightarrow PT^*(LC_0) $ and the differential map $ d\pi
:TPT^*(LC_0)\rightarrow TLC_0 $ of $\pi .$ For any $X\in
TPT^*(LC_0),$ there exists an element $\alpha\in T^*(LC_0)$ such
that $\tau (X)=[\alpha ].$  For an element $V\in T_x(LC_0),$ the
property $\alpha (V)=0$ does not depend on the choice of
representative of the class $[\alpha ].$ Thus we can define the
canonical contact structure on $PT^*(LC_0)$ by
$$
K=\{X\in TPT^*(LC_0)|\tau (X)(d\pi (X))=0\}.
$$
\par
Since we consider the coordinate $(v_2,v_3,v_4),$ we have the
trivialization $ PT^*(LC_0)\cong LC_0\times P({\Bbb R}^2)^*, $ we
call
$$
((v_2,v_3,v_4),[\xi _2:\xi _3:\xi _4])
$$
{\it a homogeneous coordinate}, where $ [\xi _2:\xi _3:\xi _4] $ is
the homogeneous coordinate of the dual projective space $P({\Bbb
R}^2)^*.$
\par
It is easy to show that $X\in K_{(x,[\xi])}$ if and only if $
\sum_{i=2}^4 \mu _i\xi _i=0, $ where $ d\pi (X)=\sum_{i=2}^4 \mu
_i\frac{\partial}{\partial v_i}. $ An immersion  $i:L\rightarrow
PT^*(LC_0)$ is said to be {\it a Legendrian immersion} if
$\text{dim}\, L=2$ and $di_q(T_qL)\subset K_{i(q)}$ for any $q\in
L.$ We also call the map $\pi\circ i$ {\it the Legendrian map} and
the set $W(i)=\text{image}\, \pi\circ i$ {\it the wave front} of
$i.$ Moreover, $i$ (or, the image of $i$) is called {\it the
Legendrian lift } of $W(i).$
\par
In order to study the lightcone pedal surface, we give a quick
survey on the Legendrian singularity theory mainly due to
Arnol'd-Zakalyukin \cite{Arnold1,Zak}. Although the general theory
has been described for general dimension, we only consider the
$3$-dimensional case for the purpose. Let $F:({\mathbb
R}^k\times{\mathbb R}^3,\bo )\longrightarrow ({\mathbb R},\bo )$ be
a function germ. We say that $F$ is {\it a Morse family} if the
mapping
$$
\Delta^*F=\left(F,\frac{\partial F}{\partial q_1}, \dots
,\frac{\partial F}{\partial q_k}\right):({\mathbb R}^k\times
{\mathbb R}^3,\bo )\lon ({\mathbb R}\times {\mathbb R}^k,\bo )
$$
is non-singular, where $(q,x)=(q_1,\dots ,q_k,x_1,x_3 ,x_3)\in
({\mathbb R}^k\times{\mathbb R}^3,\bo ).$ In this case we have a
smooth $2$-dimensional submanifold
$$
\Sigma _*(F)=\left\{ (q,x)\in ({\mathbb R}^k\times{\mathbb R}^3,\b0
)\  |\ F(q,x)=\frac{\partial F}{\partial q_1}(q,x)=\cdots
=\frac{\partial F}{\partial q_k}(q,x)=0\ \right\}
$$
and the map germ $\Phi _F:(\Sigma _*(F), \bo)\lon PT^*{\mathbb R}^3$
defined by
$$
\Phi _F(q,x)=\left(x,[\frac{\partial F}{\partial
x_1}(q,x):\frac{\partial F}{\partial x_2}(q,x): \frac{\partial
F}{\partial x_3}(q,x)]\right)
$$
is a Legendrian immersion. Then we have the following fundamental theorem of Arnol'd-Zakalyukin \cite{Arnold1,Zak}.

\begin{Pro}\rm{
All Legendrian submanifold germs in $PT^*{\mathbb R}^3$ are
constructed by the above method.}
\end{Pro}
\par
We call $F$ {\it a generating family} of $\Phi _F.$
Therefore the corresponding wave front is
$$
W(\Phi _F)\! =\!\left\{x\in {\mathbb R}^3\ |{\rm there\ exists}\
q\in {\mathbb R}^k\ {\rm such\ that}\ F(q,x)=\frac{\partial
F}{\partial q_1}(q,x)=\cdots =\frac{\partial F}{\partial
q_k}(q,x)=0\ \right\}.
$$
By definition, we have ${\cal D}_F=W(\Phi _F).$ By the previous
arguments, the lightcone pedal surface  $LP^{\pm}_M$ is the
discriminant set of the extended Lorentzian lightcone height
function $\widetilde{H}.$ We have the following proposition.

\begin{Pro}\rm{
The extended Lorentzian lightcone height function $\widetilde{H}$ is
a Morse family.}
\end{Pro}
\demo
We define another family of function
$$
\bar{H}:U\times S_t^1\times S_s^2\times \R\lon \R
$$
by
$\bar{H}((x,y),\bw ,r)=\langle \bX (x,y),\bw\rangle -r.$
We consider a $C^\infty$-diffeomorphism
$$
\Phi :U\times S_t^1\times S_s^2\times \R\lon LC_0
$$
defined by $\Phi ((x,y),\bw ,r)=((x,y),r\bw ).$ Then we have
$\widetilde{H}=\bar{H}\circ \Phi .$ It is enough to show that
$\bar{H}$ is a Morse family. For any $\bw =(\cos\theta,\sin\theta,w
_3,w_4 )\in S_t^1\times S_s^2,$ we have $w_3=\sqrt{1-w_4^2},$ so
that
$$
\bar{H}((x,y) ,\bw,r)=-x_1(p)\cos\theta -x_2(p
)\sin\theta+x_3(p)\sqrt{1-w_4^2}+x_4(p)w_4-r,
$$
where $\bX(x,y )=\bX(p )=(x_1(p),x_2(p),x_3(p), x_4(p)).$ We now
prove that the mapping
$$
\Delta^*\bar{H}=\left(\bar{H}, \frac{\partial \bar{H}}{\partial
x},\frac{\partial \bar{H}}{\partial y} \right)
$$
is non-singular at $w\in {\mathcal D}_{\bar{H}} $.
The Jacobian matrix of $\Delta ^*\bar{H}$ is given as follows:

\begin{eqnarray*}
\left(
\begin{array}{ccccc}
\langle \bX _x,\bw\rangle &  \langle \bX _y,\bw\rangle &
x_{1}\sin\theta -x_{2}\cos\theta  & {-x_3\frac{w_4}{w_3}+x_4} &
-1 \\
\langle \bX _{xx},\bw\rangle &  \langle \bX _{xy},\bw\rangle &
x_{1,x}\sin\theta -x_{2,x}\cos\theta &
{-x_{3,x}\frac{w_4}{w_3}+x_{4,x}} &
0 \\
\langle \bX _{xy},\bw\rangle &  \langle \bX _{yy},\bw\rangle&
x_{1,y}\sin\theta -x_{2,y}\cos\theta  &
{-x_{3,y}\frac{w_4}{w_3}+x_{4,y}} &
0 \\
\end{array}
\right) .
\end{eqnarray*}
By a straight forward calculation,  the determinant of the matrix
\begin{eqnarray*}
A=\left(
\begin{array}{cc}
x_{1,x}\sin\theta -x_{2,x}\cos\theta  &
\displaystyle{-x_{3,x}\frac{w_4}{w_3}+x_{4,x}} \\
x_{1,y}\sin\theta -x_{2,y}\cos\theta  &
\displaystyle{-x_{3,y}\frac{w_4}{w_3}+x_{4,y}}\\
\end{array}
\right)
\end{eqnarray*}
is equal to
$$
\frac{1}{2w_3}\vmatrix
\sin\theta & \cos\theta&-w_3 &w_4\\
\\sin\theta &-\cos\theta&-w_3 &w_4\\
x_{1,x}&x_{2,x}&x_{3,x}&x_{4,x}\\
x_{1,y}&x_{2,y}&x_{3,y}&x_{4,y}\\
\endvmatrix .$$
If det$A= 0$ then $(\cos\theta,\sin\theta,-w_3,w_4)\in T_pM$. So
that $(\cos\theta,\sin\theta,-w_1,w_2)\in N_pM\cap S_t^1\times
S_s^2$. It is impossible because
$w=(\cos\theta,\sin\theta,w_1,w_2)\in N_pM\cap S_t^1\times S_s^2$
and $N_pM$ is a Lorentzian 2-plane. Hence det$A\not=0.$ \enD
\par

By proposition 2.3, we remark that the lightcone pedal surface
$LP^\pm_M$ are wave fronts and the extended Lorentzian lightcone
height function $\widetilde{H}$ gives generating families of the
Legendrian lifts of $LP^\pm_M.$

\section{Contact with lightlike hyperplanes}
\par
In this section we consider the geometric meanings of the
singularities of the $S_t^1\times S_s^2$-valued lightcone Gauss map
(respectively, the lightcone pedal surface ) of $\bX(U)=M.$ We
consider the contact between Lorentzian surface  and lightlike
hyperplane like as the classical differential geometry. In the first
place, we briefly review the theory of contact due to Montaldi [20].
Let $X_i, Y_i$ $(i=1,2)$ be submanifolds of ${\Bbb R}^n$ with
$\text{dim}\, X_1=\text{dim}\, X_2$ and $\text{dim}\,
Y_1=\text{dim}\, Y_2.$ We say that {\it the contact of} $X_1$ and
$Y_1$ at $y_1$ is same type as {\it the contact of} $X_2$ and $Y_2$
at $y_2$ if there is a diffeomorphism germ $\Phi :({\Bbb
R}^n,y_1)\longrightarrow ({\Bbb R}^n,y_2)$ such that $\Phi
(X_1)=X_2$ and $\Phi (Y_1)=Y_2.$ In this case we write $
K(X_1,Y_1;y_1)=K(X_2,Y_2;y_2). $ It is clear that in the definition
${\Bbb R}^n$ could be replaced by any manifold. In his paper [20],
Montaldi gives a characterization of the notion of contact by using
the terminology of singularity theory.
\begin{Th} Let $X_i, Y_i$ $(i=1,2)$ be submanifolds of ${\Bbb R}^n$ with $\text{dim}\, X_1=\text{dim}\, X_2$
and $\text{dim}\, Y_1=\text{dim}\, Y_2.$
Let $g_i:(X_i,x_i)\longrightarrow
 ({\Bbb R}^n, y_i)$ be immersion germs and $f_i:({\Bbb R}^n,y_i)\longrightarrow ({\Bbb R}^p,0)$
be submersion germs with $(Y_i,y_i)=(f_i^{-1}(0),y_i).$
Then
$$K(X_1,Y_1;y_1)=K(X_2,Y_2;y_2)$$ if and only if $ f_1\circ g_1$ and $f_2\circ g_2$ are ${\mathcal K}$-equivalent.
\end{Th}

\par
We now consider a function $ {\mathcal H}:\Bbb R^4_2\times
LC_0\longrightarrow {\Bbb R} $ defined by ${\mathcal
H}(\bx,\bv)=\langle\bx ,\widetilde{\bv}\rangle -\sqrt{v_1^2+v_2^2}.$
For any $\bv _0\in LC_0,$ we denote that
$\frak{h}_{v_0}(\bx)={\mathcal H}(\bx ,\bv_0)$ and we have a
lightlike hyperplane $\frak{h}_{v_0}^{-1}(0)=
LHP(\widetilde{\bv}_0,\sqrt{v_{0,1}^2+v_{0,2}^2}).$ For any
$p_0=(x_0,y_0)\in U,$ we consider the lightlike vector $\bv
^\pm_0=\be_1\pm \be_2 (p_0)$ and $c^\pm =\langle \bX
(p_0),\widetilde{\bv_0} ^\pm \rangle ,$ then we have
$$
\frak{h}_{v^\pm_0}\circ\bX (p_0)={\mathcal H}\circ (\bX\times
id_{LC_0})(p_0,\bv^\pm_0) =H((x_0,y_0),\widetilde{\bv
_0}^\pm)-c^\pm=0.
$$
We also have relations that
$$
\frac{\partial \frak{h}_{v^\pm_0}\circ\bX}{\partial x}(p_0)=
\frac{\partial H}{\partial x}(p_0,\bv_0 ^\pm )=0,
$$
and
$$ \frac{\partial \frak{h}_{v^\pm_0}\circ\bX}{\partial y}(p_0)=
\frac{\partial H}{\partial y}(p_0,\bv_0 ^\pm )=0
$$
This means that the lightlike hyperplane
$\frak{h}_{v^\pm_0}^{-1}(0)=LHP(\widetilde{\bv} ^\pm_0,c^\pm )$ is
tangent to $M=\bX (U)$ at $p_0.$ In this case, we call each
$LHP(\widetilde{\bv} ^\pm_0,c^\pm )$ the {\it tangent lightlike
hyperplane} of $M=\bX(U)$ at $p_0=\bX(x_0,y_0).$ Moreover, the
intersection
$$
LHP(\widetilde{\bv} ^+_0,c^+)\cap LHP(\widetilde{\bv} ^-_0,c^- )
$$
is the tangent plane of $M$ at $p_0.$
Let $\bv _1,\bv _2$ be lightlike vectors. If $\bv _1,\bv _2$ are linearly dependent, then
corresponding lightlike hyperplanes $LHP(\bv _1,c_1)$ and $LHP(\bv _2,c_2)$ are parallel.
Then we have the following simple lemma.
\begin{Lem}\rm{
Let $\bX :U\longrightarrow \Bbb R^4_2$ be an immersion with $\bX
(U)$ is a
 Lorentzian  surface  and $\sigma =\pm .$
Consider two points $p_1=\bX(x_1,y_1),p_2=\bX(x_2,y_2).$ Then we
have the following assertions:
\par\noindent
{\rm (1)} $LG^\sigma _M (p_1)=LG^\sigma _M (p_2)$ if and only if $LHP(\bv _1^\sigma ,c_1^\sigma )$ and $LHP(\bv_2^\sigma ,c_2^\sigma ) $
are parallel.
\par\noindent
{\rm (2)} $LP^\sigma _M (p_1)=LP^\sigma _M (p_2)$ if and only if $LHP(\bv _1^\sigma ,c_1^\sigma )=LHP(\bv_2^\sigma ,c_2^\sigma ).$
\par\noindent
Here, $v^\pm_i=\widetilde{\be_1\pm \be_2 }(p_i)$ and $c_i^\pm
=\langle \bX (x_i,y_i),\bv _i^\pm\rangle $ for $i=1,2.$}
\end{Lem}

\par
On the other hand, for any map $f:N\lon P,$ we denote $\Sigma (f)$
the set of singular points of $f$ and $D(f)=f(\Sigma (f)).$ In this
case we call $f|_{\Sigma (f)}:\Sigma (f)\lon D(f)$ {\it the critical
part of } the mapping $f.$ For any Morse family $F:({\mathbb
R}^k\times{\mathbb R}^3,\bo )\longrightarrow ({\mathbb R},\bo ),$
$(F^{-1}(0),\bo )$ is a smooth hypersurface, so that we define a
smooth map germ $\pi _F :(F^{-1}(0),\bo )\lon ({\mathbb R}^3,0)$ by
$\pi _F(q,x)=x.$ We can easily show that $\Sigma _*(F)=\Sigma (\pi
_F).$ Therefore, the corresponding Legendrian map $\pi\circ \Phi _F$
is the critical part of $\pi _F.$
\par
We now introduce an equivalence relation among Legendrian immersion germs.
 Let  $i: (L,p) \subset (PT^*{\mathbb R}^3,p)$  and
 $i' : (L',p') \subset (PT^*{\mathbb R}^3, p')$  be Legendrian immersion germs. Then we say that  $i$  and  $i'$ are
{\it Legendrian equivalent} if there exists a contact diffeomorphism
germ $H : (PT^*{\mathbb R}^3,p) \lon  (PT^*{\mathbb R}^3,p')$  such
that  $H$  preserves fibers of  $\pi$ and that $H(L) = L'$. A
Legendrian immersion germ into $PT^*{\mathbb R}^3$ at a point is
said to be Legendrian stable if for every map with the given germ
there is a neighbourhood in the space of Legendrian immersions (in
the Whitney $C^{\infty}$ topology) and a neighbourhood of the
original point such that each Legendrian immersion belonging to the
first neighbourhood has in the second neighbourhood a point at which
its germ is Legendrian equivalent to the original germ.
\par
Since the Legendrian lift $i: (L,p) \subset (PT^*{\mathbb R}^3,p)$
is uniquely determined on the regular part of the wave front $W(i),$
we have the following simple but significant property of Legendrian
immersion germs:
\begin{Pro}\rm{
Let  $i: (L,p) \subset (PT^*{\mathbb R}^3,p)$  and
 $i': (L',p') \subset (PT^*{\mathbb R}^3, p')$  be Legendrian immersion germs such that regular sets of $\pi\circ i, \pi\circ i'$
 are dense respectively.
 Then $i, i'$ are Legendrian equivalent if and only if wave
 front sets $W(i), W(i')$ are diffeomorphic as set germs.}
 \end{Pro}
This result has been firstly pointed out by Zakalyukin [27]. The
assumption in the above proposition is a generic condition for
$i,i'.$ Especially, if $i,i'$ are Legendrian stable, then these
satisfy the assumption.
\par
We can interpret the Legendrian equivalence  by using the notion of
generating families. We denote ${\cal E}_n$ the local ring of
function germs $({\mathbb R}^n,\bo )\lon {\mathbb R}$ with the
unique maximal ideal $\mathfrak{M}_n=\{h\in {\cal E}_n\ |\ h(0)=0\
\}.$
\par
Let  $F,G : ({\mathbb R}^k\times {\mathbb R}^n,\bo ) \lon ({\mathbb
R},\bo )$  be function germs. We say that  $F$  and  $G$ are $
P$-${\cal K}$-{\it equivalent} if there exists a diffeomorphism germ
$\Psi  : ({\mathbb R}^k\times {\mathbb R}^n,\bo ) \longrightarrow
({\mathbb R}^k\times {\mathbb R}^n,\bo )$ of the form  \\
$\Psi (x,u)
= (\psi _1(q,x), \psi _{2}(x))$  for $(q,x) \in ({\mathbb R}^k\times
{\mathbb R}^n,\bo )$  such that $\Psi ^{*}(\langle F\rangle_{ {\cal
E}_{k+n}}) = \langle G\rangle_{{ \cal E}_{k+n}}$. Here  $\Psi ^{*} :
{\cal E}_{k+n} \lon {\cal E}_{k+n}$ is the pull back ${\mathbb
R}$-algebra isomorphism defined by $\Psi ^{*}(h) =  h\circ \Psi$ .
\par
Let $F: ({\mathbb R}^k\times {\mathbb R}^3,\bo )\longrightarrow
({\mathbb R},\bo )$ a function germ. We say that $F$ is a ${\cal
K}$-{\it versal deformation of} $f = F\mid_{ {\mathbb R}^k\times\{ 0
\}}$ if
$${\cal E}_k =
T_e({\cal K})(f) + \left\langle \frac{\partial F}{\partial
x_1}\mid_{ {\mathbb R}^k\times\{ \bo \}}, \frac{\partial F}{\partial
x_2}\mid_{ {\mathbb R}^k\times\{ \bo \}},\frac{\partial F}{\partial
x_3}\mid_{ {\mathbb R}^k\times\{ \bo \}}
 \right\rangle_{\mathbb R}
$$
where
$$
T_e({\cal K})(f) = \left\langle \frac{\partial f}{\partial q_1},
\dots, \frac{\partial f}{\partial q_k}, f\right\rangle_{{\cal
E}_k}.$$ (See [9].)
\par
The main result in Arnol'd-Zakalyukin's theory\cite{Arnold1,Zak} is as follows:

\begin{Th}\rm{ Let
$F, G : ({\mathbb R}^k\times {\mathbb R}^3,\bo ) \lon  ({\mathbb
R},0)$  be Morse families. Then
\par\noindent
 {\rm (1)} $\Phi _{F}$  and  $\Phi _{G}$  are Legendrian equivalent if and only if
 $F,\ G$  are  $P$-${\cal K}$-equivalent,
 \par
\par\noindent
{\rm (2)} $\Phi _{F}$  is Legendrian stable if and only if  $F$  is
a ${\cal K}$-versal deformation of  $F\mid_{ {\mathbb R}^k\times
\{\bo \}}.$}
\end{Th}
\par
Since $F,G$ are function germs on the common space germ $({\mathbb
R}^k \times {\mathbb R}^3,\bo ),$ we do no need the notion of stably
$P$-${\cal K}$-equivalences under this situation . By the uniqueness
result of the ${\cal K}$-versal deformation of a function germ,
Lemma 3.2 and Proposition 3.3, we have the following classification
result of Legendrian stable germs. For any map germ $f:({\mathbb
R}^n,\bo)\lon ({\mathbb R}^p,\bo),$ we define {\it the local ring
of} $f$ by $Q(f)={\cal E}_n/f^*(\mathfrak{M}_p){\cal E}_n.$
\begin{Pro}\rm{{\small$^{\cite{Izu2}}$}
 Let $F, G : ({\mathbb
R}^k\times {\mathbb R}^3,\bo ) \lon  ({\mathbb R},0)$  be Morse
families. Suppose that $\Phi _F,\Phi _G$ are Legendrian stable. Then
the following conditions are equivalent.
\par
{\rm (1)} $(W(\Phi _F),\bo)$ and $(W(\Phi _G),\bo )$ are diffeomorphic as germs.
\par
{\rm (2)} $\Phi _F$ and $\Phi _G$ are Legendrian equivalent.
\par
{\rm (3)} $Q(f)$ and $Q(g)$ are isomorphic as ${\mathbb
R}$-algebras, where $f=F|_{{\mathbb R}^k\times \{\bo\}},\
g=G|_{{\mathbb R}^k\times\{\bo\}}.$}
\end{Pro}

We now have tools for the study of the contact between Lorentzian
surface  and lightlike hyperplanes. Let $LP^\sigma _{M,i} :(U,
(x_i,y_i))\lon (LC_0,\bv ^\sigma  _i)$ $(i=1,2)$ be two
 lightcone pedal surface  germs of Lorentzian surface  germs $\bX _i:(U,(x_i,y_i))\lon (\R^4_2,p_i),$
where $\sigma  =\pm .$ We say that $LP^\sigma  _{M,1}$ and
$LP^\sigma  _{M,2}$ are {\it ${\cal A}$-equivalent\/} if there exist
diffeomorphism germs $\phi :(U,(x_1,y_1))\lon (U,x_2,y_2))$ and
$\Phi :(LC_0,\bv ^\sigma  _1)\lon (LC_0,\bv ^\sigma  _2)$ such that
$\Phi\circ LP^\sigma  _{M,1}=LP^\sigma  _{M,2}\circ\phi.$ If the
both of the regular sets of $LP^\sigma  _{M,i}$ are dense in
$(U,(x_i,y_i)),$ it follows from  proposition 3.5 that $LP^\sigma
_{M,1}$ and $LP^\sigma  _{M,2}$ are ${\cal A}$-equivalent if and
only if the corresponding Legendrian lift germs  are Legendrian
equivalent. This condition is also equivalent to the condition that
two generating families $\widetilde{H}_1$ and $\widetilde{H}_2$ are
$P$-${\cal K}$-equivalent by Theorem 3.4. Here,
$\widetilde{H}_i:(U\times LC_0,((x_i,y_i),\bv ^\sigma _i))\lon
{\mathbb R}$ is the extended Lorentzian lightcone height function
germ of $\bX _i.$
\par
On the other hand, we denote that $\widetilde{h}_{i,v ^\sigma
_i}(u)=\widetilde{H}_i(u,\bv ^\sigma  _i),$ then we have
$\widetilde{h}_{i,v^\pm_i}(u)=\mathfrak{h}_{v^\pm_i}\circ\bX _i(u).$
By Theorem 3.1, $K(\bX _1(U),LHP(\bv_1 ^\sigma  ,-1),\bv^\sigma
_1)=K(\bX _2(U),LHP(\bv_2^\sigma   ,-1),\bv^\sigma  _2)$ if and only
if $\widetilde{h}_{1,v_1}$ and $\widetilde{h}_{1,v_2}$ are
${\mathcal K}$-equivalent. Therefore, we can apply the previous
arguments to our situation. We denote $Q^\sigma   (\bX ,(x_0,y_0))$
the local ring of the function germ $\widetilde{h}_{v^\sigma  _0}:
(U,(x_0,y_0))\lon {\mathbb R},$ where $\bv^\sigma  _0=LP_M^\sigma
(x_0,y_0).$ We remark that we can explicitly write the local ring as
follows:
$$
Q^\pm (\bX ,(x_0,y_0))= \frac{C^\infty
_{(x_0,y_0)}(U)}{\displaystyle{\langle \langle \bX (x,y),
\widetilde{\be _1\pm \be _2}(x_0,y_0)\rangle -1\rangle _{C^\infty
_{(x_0,y_0)}(U)}}},
$$
where $C^\infty _{(x_0,y_0)}(U)$ is the local ring of function germs
at $(x_0,y_0)$ with the unique maximal ideal
$\mathfrak{M}_{(x_0,y_0)}(U).$
\begin{Th}\rm{
Let $\bX _i:(U,(x_i,y_i)) \lon (\R^4_2,\bX _i(x_i,y_i))$ $(i=1,2)$
be an immersion germs with $\bX(U)=M$ is a Lorentzian surface such
that the corresponding Legendrian lift germs are Legendrian stable
and $\sigma =\pm .$ Then the following conditions are equivalent:
\par
{\rm (1)} lightcone pedal surface  germs $LP_{M,1}^\sigma  $ and
$LP_{M,2}^\sigma $ are ${\cal A}$-equivalent.
\par
{\rm (2)} $\widetilde{H}_1$ and $\widetilde{H}_2$ are $P$-${\cal K}$-equivalent.
\par
{\rm (3)} $\widetilde{h}_{1,v_1}$ and $\widetilde{h}_{1,v_2}$ are ${\cal K}$-equivalent.
\par
{\rm (4)} $K(\bX _1(U),LHP(\bv^\sigma
_1,c_1^\sigma),\bv^{\sigma}_{1})=K(\bx _2(U), LHP(\bv^\sigma
_2,c_2^\sigma),\bv^{\sigma}_{2})$
\par
{\rm (5)} $Q^\sigma  (\bX _1 ,(x_1,y_1))$ and $Q^\sigma  (\bX _2
,(x_2,y_2))$ are isomorphic as ${\mathbb R}$-algebras.}
\end{Th}
\demo
By the previous arguments (mainly by Theorem 3.1),
it has been already shown that conditions (3) and (4) are equivalent.
Other assertions follow from proposition 3.5.
\enD

\par
Given an immersion germ $\bX :(U,(x_0,y_0))\lon (\R^4_2,\bX
(x_0,y_0))$ with $\bX(U)=M$ is a Lorentzian surface, we call each
set
$$
(\bX ^{-1}(LHP(\bv^\pm ,c^\pm )),(x_0,y_0))
$$
a {\it tangent lightlike hyperplane indicatrix germ} of $\bX ,$
where $\bv ^\pm =\be _1\pm\be _2(x_0,y_0)$ and $c^\pm=\langle
\bX(x_0,y_0),\bv ^\pm \rangle.$ Moreover, by the above results, we
can borrow some basic invariants from the singularity theory on
function germs. We need ${\cal K}$-invariants for function germ. The
local ring of a function germ is a complete ${\cal K}$-invariant for
generic function germs. It is, however, not a numerical invariant.
The ${\cal K}$-codimension (or, Tyurina number) of a function germ
is a numerical ${\cal K}$-invariant of function germs
\cite{martine}. We denote that
$$
\mbox{\rm L-ord}^\pm (\bX ,(x_0,y_0))={\rm dim}_\R\, \frac{C^\infty
_{(x_0,y_0)}(U)}{\langle
\widetilde{h}_{v_0^\pm}(x,y),\widetilde{h}_{v_0^\pm ,x}(x,y),
\widetilde{h}_{v_0^\pm ,y}(x,y)\rangle }.
$$
Usually $\mbox{\rm L-ord}^\sigma  (\bx ,u_0)$ is called {\it the
${\cal K}$-codimension of} $\widetilde{h}_{v^\sigma _0},$ where
$\sigma =\pm .$ However, we call it the {\it order of contact with
the tangent lightlike hyperplane} at $\bX (x_0,y_0).$ We also have
the notion of corank of function germs.
$$
\mbox{\rm L-corank}^\sigma  (\bX ,(x_0,y_0))= 2-{\rm rank}\, {\rm
Hess}(\widetilde{h}_{v^\sigma _0}(x_0,y_0)),
$$
where $\bv^\pm _0=\be _1\pm\be _2 (x_0,y_0).$
\par
By proposition 2.1, $\bX (x_0,y_0)$ is a $L^\sigma $-parabolic point
if and only if $\mbox{\rm L-corank}^\sigma  (\bX ,(x_0,y_0)) \geq
1.$ Moreover $\bX (x_0,y_0)$ is a lightlike umbilic point if and
only if $\mbox{\rm L-corank}^\sigma  (\bX ,(x_0,y_0))=2.$
\par
On the other hand, a function germ $f:({\mathbb R}^{n-1},\ba )\lon
{\mathbb R}$ has the $A_k$-type singularity if $f$ is ${\cal
K}$-equivalent to the germ $\pm u_1^2\pm \cdots \pm
u_{n-2}^2+u_{n-1}^{k+1}.$ If $\mbox{\rm L-corank}^\sigma  (\bX
,(x_0,y_0))=1,$ the extended Lorentzian lightcone height function
$\widetilde{h}_{v^\sigma  _0}$ has the $A_k$-type singularity at
$(x_0,y_0)$ in generic. In this case we have $\mbox{\rm
L-ord}^\sigma  (\bx ,u_0)=k.$ This number $k$ is equal to the order
of contact in the classical sense (cf., \cite{Bruce-Giblin1}). This
is the reason why we call $\mbox{\rm L-ord}^\sigma  (\bX
,(x_0,y_0))$ the order of contact with the tangent lightlike
hyperplane at $\bX (x_0,y_0).$
\par
As a corollary of the theorem 3.6,
we have the following result.
\begin{Co}\rm{
Under the assumptions of Corollary 3.6, if the tangent lightlike
hyperplane indicatrix germ $LHP_{M_1}^\sigma  $ and
$LHP_{M_2}^\sigma $ are ${\cal A}$-equivalent, then tangent
lightcone indicatrix germs
$$(\bX ^{-1}(LHP(\bv_1^\pm ,c_1^\pm )),(x_1,y_1))\quad {\rm and }\quad
(\bX ^{-1}(LHP(\bv_2^\pm ,c_2^\pm )),(x_2,y_2))
$$ are diffeomorphic as set germs.}
\end{Co}
\demo
Notice that the tangent lightlike hyperplane indicatrix germ
 of $\bX_i$ is the zero
level set of $h_{i,\lambda _i}.$ Since ${\mathcal K}$-equivalence
among function germs preserves the zero-level sets of function
germs, the assertion follows from theorem 3.6. \enD
\section{Classification of singularities of $S_t^1\times S_s^2$-valued
lightcone Gauss maps and lightcone pedalsurfaces
}
\par In this section we consider generic singularities of
$S_t^1\times S_s^2$-valued lightcone Gauss maps and lightcone pedal surfaces.
We consider the space of Lorentzian embeddings ${\rm Emb}_L\,
(U,\R^4_2)$ with Whitney $C^\infty$-topology, where $U\subset \R^2$
is an open subset. We  have the following theorem.
\begin{Th}\rm{
There exists an open dense subset ${\cal O}\subset {\rm
Emb}_L\,(U,\R^4_2)$ such that for any $\bX\in {\cal O},$ the
following conditions hold:
\par
{\rm (1)} Each  lightlike parabolic set ${\mathcal K}_l(1,\sigma
1)^{-1}(0)$ is a regular curve. We call such a curve {\it the
lightlike parabolic curve}.
\par
{\rm (2)} The lightcone pedal surface $LP_M^\sigma  $ along the
lightlike parabolic curve is the cuspdialedge except isolated
points. At these points $LP_M^\sigma  $ is the swallowtail.
\par
Here, $\sigma  =\pm$ and a map germ $f:({\mathbb R}^2,\ba )\lon
({\mathbb R}^3,\bb)$ is called {\it a cuspidaledge} if it is ${\cal
A}$-equivalent to the germ $(u_1,u_2^2,u_2^3)$ {\rm (}cf., {\rm Fig.
1)} and {\it a swallowtail} if it is ${\cal A}$-equivalent to the
germ $(3u_1^4+u_1^2u_2,4u_1^3+2u_1u_2,u_2)$ .}
\end{Th}
\par
For the proof of Theorem 4.1, we consider the function ${\mathcal
H}:\R^4_2\times LC_{0}\lon \R$ which is given in \S 3. We claim that
${\mathcal H}_v$ is a submersion for any $\bv\in LC_{0},$ where
${\mathcal H}_v(\bx)={\mathcal H}(\bx,\bv ).$ For any $\bX \in {\rm
Emb}_L\, (U,\R^4_2),$ we have $\widetilde{H}={\cal H}\circ (\bX
\times id _{LC_{0}}).$ We also have the $\ell $-jet extension
$$
j^{\ell }_1\widetilde{H}:U\times LC_{0}\lon J^{\ell }(U,{\mathbb R})
$$
defined by $j^{\ell }_1\widetilde{H}((x,y),\bv)=j^{\ell
}\widetilde{h}_v (x,y).$ We consider the trivialization
$J^{\ell}(U,{\mathbb R})\equiv U\times {\mathbb R}\times J^{\ell
}(2,1).$ For any submanifold $Q\subset J^\ell (2,1),$ we denote that
$\widetilde{Q}=U\times\{0\}\times Q.$ Then we have the following
proposition as a corollary of Lemma 6 in Wassermann
\cite{Wassermann}. (See also Montaldi \cite{Montaldi2} and Looijenga
\cite{Looijenga}).
\begin{Pro}\rm{
Let $Q$ be a submanifold of $J^{\ell }(2,1).$ Then the set
$$
T_Q=\{\bx\in {\rm Emb}_L\, (U,\R^4_2)\ |\ j^{\ell}_1H\ \mbox{\rm is
transversal to}\ \widetilde{Q}\ \}
$$
is a residual subset of ${\rm Emb}_s\, (U,\R^4_2).$ If $Q$ is a
closed subset, then $T_Q$ is open.}
\end{Pro}
If we consider ${\cal K}$-orbits in $J^\ell (2,1),$ we obtain the
proof of Theorem 4.1, so that we omit the detailed discussion. The
assertion of Theorem 4.1 can be interpreted that the Legendrian lift
of the lightcone pedal surface $LP_M^\pm $ of $\bX \in {\cal O}$ is
Legendrian stable at each point. Since the Legendrian lift is the
Legendrian covering of the Lagrangian lift of  $LG_M^\pm ,$ it has
been known that the corresponding singularities of $LG_M^\pm$ are
folds or cusps \cite{Arnold1}. Hence, we have the following
corollary.
\begin{Co}\rm{
Let ${\cal O}\subset {\rm Emb}_L\, (U,\R^4_2)$ be the same open
dense subset as in Theorem 4.1. For any $\bX\in {\cal O},$ the
followings hold:
\par
{\rm (1)} A lightlike parabolic point $(x_0,y_0)\in U$ is a fold of
the $S_t^1\times S_s^2$-valued lightcone Gauss map $LG_M^\sigma  $
if and only if it is a cuspidaledge of the lightcone pedal surface
$LP_M^\sigma  .$
\par
{\rm (2)} A lightlike parabolic point $(x_0,y_0)\in U$ is a cusp of
the $S_t^1\times S_s^2$-valued lightcone Gauss map $LG_M^\sigma  $
if and only if it is a swallowtail of the lightcone pedal surface
$LP_M^\sigma  .$
\par
Here, a map germ $f:({\mathbb R}^2,\ba)\lon ({\mathbb R}^2,\bb)$ is
called {\it a fold} if it is ${\cal A}$-equivalent to the germ
$(u_1,u_2^2)$ and {\it a cusp} if it is ${\cal A}$-equivalent to the
germ $(u_1,u_2^3+u_1u_2).$}
\end{Co}
\par
Following the terminology of Whitney \cite{whitney}, we say that a
surface $\bX :U\lon \R^4_2$ has {\it the excellent lightcone pedal
surface} $LP_M^\sigma $ if the Legendrian lift of $LP^\sigma _M$ is
a stable Legendrian immersion at each point. In this case, the
lightcone pedal surface $LP_M^\sigma $ has only cuspidaledges and
swallowtails as singularities. Proposition 4.1 asserts that a
Lorentzian surface  with the excellent lightcone pedal surface is
generic in the space of all Lorentzian surface  in $\R^4_2 .$ We now
consider the geometric meanings of cuspidaledges and swallowtails of
the lightcone pedal surface. We have the following results analogous
to the results of Banchoff et al \cite{Banchoff}.
\begin{Th}\rm{
Let $LP_M^\sigma :(U,(x_0,y_0))\lon (\R^4_2,p_0 )$ be the excellent
lightcone pedal surface  of a Lorentzian surface  $\bX$ and
$\widetilde{h}_{v^\sigma_0} :(U,(x_0,y_0))\lon \R$ be the extended
lightcone height function germ at $\bv ^\pm_0=\be_1\pm\be _2(p_0),$
where $\sigma =\pm.$ Then we have the following:
\par\noindent
{\rm (1)} $(x_0,y_0)$ is a lightlike parabolic point of $\bX$ if and
only if $\mbox{\rm L-corank}^\sigma (\bX ,(x_0,y_0))=1.$
\par\noindent
{\rm (2)} If $(x_0,y_0)$ is a lightlike parabolic point of $\bX ,$
then $\widetilde{h}_{v^\sigma _0}$ has the $A_k$-type singularity
for $k=2,3.$
\par\noindent
{\rm (3)} Suppose that $(x_0,y_0)$ is a lightlike parabolic point of
$\bX.$ Then the following conditions are equivalent:
\par
{\rm (a)} $LP_M^\sigma$ has a cuspidaledge at $(x_0,y_0)$
\par
{\rm (b)} $\widetilde{h}_{v^\sigma_0}$ has the $A_2$-type
singularity.
\par
{\rm (c)} $\mbox{\rm L-ord}^\sigma (\bX ,(x_0,y_0))=2.$
\par
{\rm (d)} The tangent lightlike hyperplane indicatrix is a ordinary
cusp, where a curve $C\subset {\mathbb R}^2$ is called {\it a
ordinary cusp} if it is diffeomorphic to the curve given by
$\{(u_1,u_2)\ |\ u_1^2-u_2^3=0\ \}.$
\par
{\rm (e)} For each $\varepsilon >0,$ there exist two distinct points
$(x_i,y_i)\in U$ $(i=1,2)$ such that
$$
\|(x_0,y_0)-(x_i,y_i)\|<\varepsilon
$$ for $i=1,2,$
both of $(x_i,y_i)$ are not lightlike parabolic points and the
tangent lightlike hyperplanes to $M=\bx (U)$ at $(x_i,y_i)$ are
parallel.
\par\noindent
{\rm (4)} Suppose that $(x_0,y_0)$ is a lightlike parabolic point of
$\bX.$ Then the following conditions are equivalent:
\par
{\rm (a)} $LP_M^\sigma$ has a swallowtail at $(x_0,y_0)$
\par
{\rm (b)} $\widetilde{h}_{v^\sigma _0}$ has the $A_3$-type
singularity.
\par
{\rm (c)} $\mbox{\rm L-ord}^\sigma (\bX ,(x_0,y_0))=3.$
\par
{\rm (d)} The tangent lightlike hyperplane indicatrix is a point or
a tachnodal, where a curve $C\subset {\mathbb R}^2$ is called {\it a
tachnodal} if it is diffeomorphic to the curve given by
$\{(u_1,u_2)\ |\ u_1^2-u_2^4=0\ \}.$
\par
{\rm (e)} For each $\varepsilon >0,$ there exist three distinct
points $(x_i,y_i)\in U$ $(i=1,2,3)$ such that
$$
\|(x_0,y_0)-(x_i,y_i)\| <\varepsilon
$$
for $i=1,2,3$ and the tangent lightlike hyperplanes to $M=\bx (U)$
at $(x_i,y_i)$ are parallel.
\par
{\rm (f)} For each $\varepsilon >0,$ there exist two distinct points
$(x_i,y_i)\in U$ $(i=1,2)$ such that
$$
\|(x_0,y_0)-(x_i,y_i)\|<\varepsilon
$$
for $i=1,2$ and the tangent lightlike hyperplanes to $M=\bx (U)$ at
$(x_i,y_i)$ are equal.}
\end{Th}
\demo We have shown that $(x_0,y_0)$ is a lightlike parabolic point
if and only if
$$
\mbox{\rm L-corank}^\sigma (\bX ,(x_0,y_0))\geq 1.
$$
We have  $\mbox{\rm L-corank}^\sigma (\bX ,(x_0,y_0))\leq 2.$ Since
the extended lightcone height function germ $\widetilde{H}:(U\times
LC_0,((x_0,y_0),\bv _0))\lon {\mathbb R}$ can be considered as a
generating family of the Legendrian lift of $LP_M^\sigma ,$
$\widetilde{h}_{v^\sigma _0}$ has only the $A_k$-type singularities
$(k=1,2,3).$ This means that the corank of the Hessian matrix of
$\widetilde{h}_{v^\sigma _0}$ at a lightlike parabolic point is $1.$
The assertion (2) also follows. By the same reason, the conditions
(3);(a),(b),(c) (respectively, (4); (a),(b),(c)) are equivalent. If
the height function germ $\widetilde{h}_{v^\sigma _0}$ has the
$A_2$-type singularity, then it is ${\cal K}$-equivalent to the germ
$\pm u_1^2+u_2^3.$ Since the ${\cal K}$-equivalence preserve the
diffeomorphism type of zero level sets, the tangent lightlike
hyperplane indicatrix is diffeomorphic to the curve given by $\pm
u_1^2+u_2^3=0.$ This is the ordinary cusp. The normal form for the
$A_3$-type singularity is given by $\pm u_1^2+u_2^4,$ so that the
tangent lightlike hyperplane indicatrix is diffeomorphic to the
curve $\pm u_1^2+u_2^4=0.$ This means that the condition (3),(d)
(respectively, (4),(d)) is also equivalent to the other conditions.
\par
Suppose that $(x_0,y_0)$ is a lightlike parabolic point, then the
$S_t^1\times S_s^2$-valued lightcone Gauss map has only folds or
cusps. If the point $(x_0,y_0)$ is a fold point, there is a
neighbourhood of $(x_0,y_0)$ on which the $S_t^1\times S_s^2$-valued
lightcone Gauss map is 2 to 1 except at the lightlike parabolic
curve (i.e, fold curve). By Lemma 3.2, the condition (3), (e) is
satisfied. If the point $(x_0,y_0)$ is a cusp, the critical value
set is an ordinary cusp. By the normal form, we can understand that
the $S_t^1\times S_s^2$-valued lightcone Gauss map is 3 to 1 inside
region of the critical vale. Moreover, the point $(x_0,y_0)$ is in
the closure of the region. This means that the condition (4),(e)
holds. We can also observe that near by a cusp point, there are 2 to
1 points which approach to $(x_0,y_0).$ However, one of those points
are always lightlike parabolic points. Since other singularities do
not appear for in this case, so that the condition (3),(e)
(respectively, (4),(e)) characterizes a fold (respectively, a cusp).
\par
If we consider the lightcone pedal surface in stead of the lightcone
Gauss map, the only singularities are cuspidaledges or swallowtails.
For the swallowtail point $(x_0,y_0)$, there is a self intersection
curve approaching to $(x_0,y_0).$ On this curve, there are two
distinct point $(x_i,y_i)$ $(i=1,2)$ such that $LP_M^\sigma
(x_1,y_1)= LP_M^\sigma (x_2,y_2).$ By Lemma 3.2, this means that
tangent lightlike hyperplane to $M=\bx (U)$ at $(x_i,y_i)$ are
equal. Since there are no other singularities in this case, the
condition (4),(f) characterize a swallowtail point of $LP_M^\sigma
.$ This completes the proof. \enD

{\small
\bigskip
\par\noindent
Donghe Pei, School of Mathematics and Statistics, Northeast Normal
University, Changchun 130024, P.R.China
\par\noindent e-mail:{\tt peidh340@nenu.edu.cn}
\par\noindent
Jianguo Sun, School of Mathematics ,China University of petroleum,
Dongying 257061, P.R.China
\par\noindent
e-mail:{\tt sunjg616@163.com}
\par\noindent
Qi Wang, School of Mathematics and Statistics, Northeast Normal
University, Changchun 130024, P.R.China

\end{document}